\documentclass[12pt]{amsart}
\usepackage{graphicx}
\usepackage{amssymb}
\usepackage{epstopdf}
 \usepackage{latexsym}
\usepackage{amsmath}
\usepackage{amsfonts} 
\usepackage{amsthm}

\newtheorem{Proposition}{Proposition}
\newtheorem{Theorem}[Proposition]{Theorem}

\newtheorem{Corollary}[Proposition]{Corollary}
\newtheorem{Remark}[Proposition]{Remark}

\def\XXint#1#2#3{{\setbox0=\hbox{$#1{#2#3}{\int}$}
\vcenter{\hbox{$#2#3$}}\kern-.5\wd0}}

       \def\z{\noindent}

\def\z{\noindent}

    \def\sqr#1#2{{\vcenter{\vbox{\hrule height .#2pt
                             \hbox{\vrule width .#2pt height#1pt \kern#1pt
                                   \vrule width .#2pt}
                             \hrule height .#2pt}}}}

     \def\CC{\mathbb{C}}
    
    \def\NN{\mathbb{N}}

    \def\ZZ{\mathbb{Z}}

\def\be{\begin{equation}}
\def\ee{\end{equation}}
\def\calJ{\mathcal{J}}
\def\calB{\mathcal{B}}
\def\calH{\mathcal{H}}
\def\al{\alpha}
\def\nint{\not\negmedspace{\negthinspace\int}}
\def\ninto{\nint_{\negthickspace{0}}}

\begin{document} 

\title[Jacobi polynomials for general parameters ]{Orthogonality of Jacobi and Laguerre polynomials for general parameters via the Hadamard finite part}
\author{Rodica D. Costin}

\maketitle

\begin{abstract}

Orthogonality of the Jacobi and of Laguerre polynomials, $P_n^{(\al,\beta)}$ and $L_n^{(\al)}$, is established for $\al,\beta\in\CC\setminus\ZZ_-,\ \al+\beta\ne -2,-3,\ldots$ using the Hadamard finite part of the integral which gives their orthogonality in the classical cases. Riemann-Hilbert problems that these polynomials satisfy are found.

The results are formally similar to the ones in the classical case (when $\Re\al,\Re\beta>-1$).

\end{abstract}


\section{Introduction}

Orthogonality of classical polynomials is key for the study of many properties of these polynomials and their applications, and it is important to find constructive formulas for the bilinear functional that gives orthogonality.

Such formulas have been usually obtained by taking the analytic continuation  (in the parameters) of the Borel measure of the classical case. Carlson used an integral kernel to establish this continuation and he proves the existence of Jacobi series for general parameters \cite{Carlson}. More recently Kuijlaars, Mart\'inez-Finkelshtein and Orive find by analytic continuation that orthogonality of Jacobi polynomials can be established in some cases by integration on special paths in the complex plane; they also derive an associated Riemann-Hilbert problem \cite{KJ}; in other cases incomplete or quasi-orthogonality is found, or even multiple orthogonality conditions (see also \cite{MF-O}).

In the present paper analytic continuation is established using the Hadamard finite part of (possibly) divergent integrals. Since these can be manipulated much like integrals, the classical formulas which are analytic in the parameters are formally similar. Orthogonality of the Laguerre polynomials $L_n^{(\al)}$ is obtained in \S\ref{ortoL} for $\al\in\CC\setminus \ZZ_-$. Orthogonality of the Jacobi polynomials $P_n^{(\al,\beta)}$ is established in \S\ref{ortoJ} for $\al,\beta\in\CC\setminus \ZZ_-,\, \al+\beta\ne -2,-3,\ldots $. 

The existence of an associated Riemann-Hilbert problem for polynomials  orthogonal with respect to a Borel measure on the line, introduced in \cite{FIK}, is now a well known result and technique which has proved very useful in deducing properties of these polynomials - see \cite{Deift_book}, \cite{Ismail}, \cite{KJ_lecture}.
Associated Riemann-Hilbert problems in the present generalized context are found for the Laguerre polynomials in \S\ref{RHLP}  and  for the Jacobi polynomials in \S\ref{RHJP}.

It should be mentioned that once analyticiy in parameters of the Hadamard finite part is established (in \S\ref{osin}) the remaining results follow by analytic continuation (direct proofs are given, just as a confirmation, in the Appendix). However, classical boundary conditions of the associated Riemann-Hilbert problems are not formulated in analytic terms for non-positive parameters (see \cite{Ismail}, \cite{KJ_lecture}) and need reformulation to insure uniqueness of the solution.

\section{The Hadamard finite part of integrals $\int_0^xt^{\alpha-1}f(t)\, dt$}

The concept of the {\em{finite part}} of a (possibly divergent) integral was introduced by Hadamard \cite{Hadamard} as a convenient way to express solutions of differential equations. He showed that this finite part of an integral (which coincides with the usual value if the integral is convergent) can be combined and manipulated in much the same way as usual integrals: they are additive on the interval of integration, changes of variable are allowed, etc. (They do not behave well with respect to inequalities.) The finite part can be calculated either by Taylor series, or by integration along closed paths in the complex plane.

Subsequently the Hadamard finite part has been interpreted in terms of distributions (see, e.g. \cite{Morton-Krall}) and it turned out that many problems of mathematical physics have solutions expressible as the Hadamard finite part of (divergent) integrals, and numerical methods of calculations have been subsequently developed (see for example \cite{Davis_Rabinowitz}).

The present section contains some properties of the Hadamard finite part of integrals of the type $\int_0^xt^{\alpha-1}f(t)\, dt$ with $f$ analytic at $0$; when $\Re\alpha\leq 0$, $\alpha\not\in(-\NN)$, its Hadamard finite part is denoted here by
$$\ninto^{x}t^{\alpha-1}f(t)\, dt$$

\subsection{Notations}

$\NN=\{ 0,1,2,\ldots\}$. $\ZZ_-=\{-1,-2,-3,\ldots\}$.

\subsection{Analiticity in $\al$}\label{osin}

For $r>0$ denote by $D_0$ the disk 
$$D_0=\{x\in\CC;|x|<r\}$$
 Let $\calH(D_0)$ the Banach space of functions which are analytic on the disk $D_0$ and continuous on $\overline{D}_0$ with the sup norm.

For $\alpha\in\CC,\ \alpha\not\in-\NN$ define the operator $ \calJ_\al $ on $\calH(D_0)$ by
\be\label{defJ}
{\mbox{if\ }}f(x)=\sum_{n=0}^\infty f_nx^n\ \ {\mbox{then\ \ }} \calJ_\al f\, (x)=\sum_{n=0}^\infty \frac{f_n}{n+\al}\, x^n
\ee

\begin{Remark} $ \calJ_\al f$ solves the linear nonhomogeneous equation
\be\label{eq1e}
\frac{du}{dx}+\frac{\al}{x}\, u=\frac{1}{x}f(x)
\ee
and it is its only solution which is analytic at $x=0$.
\end{Remark}

\z This is easy to see by searching for solutions of (\ref{eq1e}) as power series at $x=0$. On the other hand, solving (\ref{eq1e}) using an integrating factor we see that

\begin{Remark}  If $\Re\al>0$ then 
\be\label{intJ}
 \calJ_\al f\, (x)=x^{-\al}\int_0^xt^{\al-1}f(t)\, dt
\ee
\end{Remark}

Proposition\,\ref{P1} shows that all the functions defined by the integral (\ref{intJ}) can be analytically continued in $\al$ except for $\al\in -\NN$ where they have poles of order one (manifestly seen in (\ref{defJ})). 


Let $0<r'<r$ and denote $D_0'=\{x\in\CC;|x|<r'\}$.

\begin{Proposition}\label{P1}

The linear operator $ \calJ_\al :\calH(D_0)\to\calH({D_0'})$ defined by (\ref{defJ}) is compact and analytic in $\al$ on $\CC\setminus(-\NN)$.

\end{Proposition}

{\em{Proof.}}

For $N\in\NN$ denote by $\calJ_\al^{[N]}$ the following finite rank and analytic in $\al$ operators on $\calH(D_0)$:
$$\calJ_\al^{[N]}f\, (x)=\sum_{n=0}^N\frac{f_n}{n+\al}\, x^n$$

Let $\delta>0$ and consider $\al$ such that dist$(\al,-\NN)\geq \delta$.

We have $|f_n|\leq \|f\| r^{-n}$, therefore if $|x|<r'$ we have
$$\big|  \calJ_\al f\, (x)-\calJ_\al^{[N]}f\, (x)\big|\leq\|f\| (r'/r)^N\sum_{k=0}^\infty \frac{(r'/r)^k}{|k+N+\al|}\leq \|f\| \frac{(r'/r)^N}{\delta\, (1-r'/r)}$$
which shows that $\calJ_\al^{[N]}$ converge to $ \calJ_\al $ in norm, and uniformly in $\al$ on compact subsets of $\CC\setminus(-\NN)$.
\qed

{\bf{Definition.}} For $f\in\calH(D_0)$ we denote
\be\label{def_int}
  \ninto^{x}\, t^{\al-1}f(t)\, dt=x^{\al} \calJ_\al f\, (x)
\ee

Note that the operator (\ref{def_int}) is a bona fide integral for $\Re \al>0$, otherwise  it is possibly  divergent. For $\Re \al\leq 0$, (\ref{def_int}) is the Hadamard finite part of the corresponding singular integral \cite{Hadamard}.

{\bf{Remark.}}  If $f(x)$ is a multiple of $x$ (or higher powers of $x$) there may be an ambiguity in the choice of $\al$ in (\ref{def_int}). For example, if $f(t)=tg(t)$ then 
\be\label{oeq1}
 \ninto^{x}\, t^{\al-1}f(t)\, dt=\ninto^x\, t^\al g(t)\, dt\, =\, x^{\al+1}\calJ_{\al+1}g\, (x)
 \ee
The two results in (\ref{def_int}) and (\ref{oeq1}) are nevertheless the same, since we have
$$\calJ_\al\left[ xg\right]=x\calJ_{\al+1}g$$

\subsection{Usual properties of integrals}
The Hadamard finite part (\ref{def_int}) satisfies usual properties of integrals: it is additive on intervals, shown in Proposition\,\ref{addi}, satisfies a (generalized)  Fundamental Theorem of Calculus, and as a consequence, integration by parts holds, shown by Proposition\,\ref{FdThCal}.

\begin{Proposition}\label{addi}

For $f\in\calH(D_0)$ and $x,\xi\in D_0'$ we have
$$\ninto^{\xi}\,  t^{\al-1} f(t)\, dt\,+\,\int_\xi^{x}\, t^{\al-1} f(t)\, dt=\ninto^{x}\,  t^{\al-1} f(t)\, dt$$
\end{Proposition}

The proof is found in \S\ref{pfaddi}.

\begin{Proposition}\label{FdThCal} 

Let $f,f',g,g'\in\calH(D_0)$ and $x\in D_0'$.

(i) We have
\be\label{fdthcal}
\ninto^{x}\, \frac{d}{dt}\left[ t^\al f(t)\right]\, dt=x^\al f(x)
\ee

(ii) As a consequence, integration by parts holds:
\be\label{ip1s}
{ \ninto}^{x}\, f(t)\, \frac{d}{dt}\left[ t^\al g(t)\right]\, dt\, =\, x^\al f(x)g(x)\, -\, { \ninto}^{x}\,  t^\al f'(t)\, g(t)\, dt 
\ee

\end{Proposition}

The proof is found in \S\ref{pfFdThCal}.

\subsection{A Riemann-Hilbert problem on half line}

Let $\al\in\CC\setminus \ZZ_-$ with (possibly) $\Re\al\leq -1$. 

If $f$ is a function analytic on\footnote{In fact, analiticity of $f$ is only required at $x=0$.} $[0,+\infty)$ and $f\in L^1(0,+\infty)$ we can define
\be\label{nintoinf}
\ninto^\infty t^\al f(t)\, dt=\ninto^\xi  t^\al f(t)\, dt\, +\, \int_\xi^\infty  t^\al f(t)\, dt
\ee
where $\xi>0$ is any positive number so that $f$ is analytic on the disk $|z|<\xi$ and the usual branch of $t^\al$ is considered (i.e. $|t^\al |=t^{\Re \al}$ for $t>0$).

{{Note}} that (\ref{nintoinf}) is analytic in $\al$. 
Note that
$$\ninto^\infty x^\al e^{-x}\, dx=\Gamma(\al+1)$$

{\bf{Notation.}} Let $\mathcal{O}_0$ be the set of germs of analytic functions at $0$. We write $f(z)=z^\al\mathcal{O}_0+{O}(1)$ for $z\to 0$ to mean that there is $f_0\in\mathcal{O}_0$ so that $f(z)-z^\al f_0(z)= {O}(1)$ for $z\to 0$. (Obviously, it can be equivalently required that $ f_0(z)$ be a polynomial of degree $\geq -\Re\al$.)

 Let $p(x)$ be a polynomial.  
 Consider the following Riemann-Hilbert problem on $(0,+\infty)$:

\z {{(i)}} $f$ is analytic on $\CC\setminus [0,+\infty)$ and the following limits exist for $x>0$:
$$f^+(x)=\lim_{z\to x;\Re z>0} f(z),\ \ \ \ \ f^-(x)=\lim_{z\to x;\Re z<0} f(z)$$

\z {{(ii)}} $\displaystyle{f^+(x)=f^-(x)+ p(x)x^\al e^{-x}}$ for $x>0$.

\z {{(iii$_\infty$)}} $\displaystyle{f(z)={O}\left(z^{-1}\right)}$ for $z\to\infty$;

\z {{(iii$_{0}$)}} $f(z)=  z^\al\mathcal{O}_0+{O}(1)$ for $z\to 0$.

\begin{Theorem}\label{ThRHsl}
\z The Riemann-Hilbert problem {\em{(i)-(iii)}} has the unique solution $f=\mathcal{C}_\al[p]$ where
\be\label{notC}
\mathcal{C}_\al[p]\, (z)=\frac{1}{2\pi i}\ninto^{\infty}\, \frac{t^\al p(t)e^{-t}}{t-z}\, dt
\ee

\end{Theorem}

{\em{Proof of {Theorem}\,\ref{ThRHsl}.}}

{\em{Existence.}} We first establish that the function (\ref{notC}) satisfies {{(i)}}-{{(iii$_\infty$)}}.

Let $z_0\in \CC\setminus [0,+\infty)$ and $x>0$. Let $\xi$ so that $0<\xi<x$ and $\xi<|z_0|$. 

Write the Hadamard finite part in (\ref{notC}) as
\be\label{descomp}
\ninto^{\infty}\, \frac{t^\al p(t)e^{-t}}{t-z}\, dt=\ninto^{\xi}\, \frac{t^\al p(t)e^{-t}}{t-z}\, dt+\int_\xi^{\infty}\, \frac{t^\al p(t)e^{-t}}{t-z}\, dt
\ee

The first term on the right side of (\ref{descomp}) is analytic in $z$ for $|z|>\xi$ (therefore it is analytic at $z_0$): we can expand $(t-z)^{-1}=-z^{-1}\sum_nt^n/z^n$, an absolutely convergent series which is obviously ${{O}\left(z^{-1}\right)}$ for $z\to\infty$.

The last integral in (\ref{descomp}) satisfies {{(i)}}-{{(iii$_\infty$)}} for $|z|>\xi $ by the classical theory (see, e.g. \cite{Ismail}, \cite{KJ_lecture}).

To establish that (\ref{notC}) satisfies {{(iii$_{0}$)}} note that  $\mathcal{C}_\beta[1]$ (for any $\beta\not\in \ZZ_-$) satisfies the differential equation
\be\label{difxbeta}
\frac{d\phi}{dz}=\left(\frac{\beta}{z}-1\right)\phi-\frac{a}{z}\ \ \ \ {\mbox{where\ }} a=\frac{\Gamma(\beta+1)}{2\pi i}
\ee
 (see \S\ref{pfdifeq} for details) having the general solution $\phi(z)=\phi_{an}(z)+Cz^\beta e^{-z}$ with $\phi_{an}$ a function analytic at $z=0$. Therefore $\mathcal{C}_\beta[1]= z^\beta\mathcal{O}_0+{O}(1)$ for $z\to 0$. If $p(x)$ is a polynomial with $p(x)=\sum_kp_kx^k$ then $\mathcal{C}_\al[p]=\sum_kp_k\mathcal{C}_{\al+k}[1]= z^\al\mathcal{O}_0+{O}(1)$ for $z\to 0$.

{\em{Uniqueness.}} Let $f$ be any function satisfying {{(i)}}-{{(iii)}}. Then $g=f-\mathcal{C}_\al[p]$
is analytic on $\CC\setminus [0,+\infty)$ and $g^+(x)=g^-(x)$ for $x>0$, therefore $g$ is analytic on $\CC\setminus \{0\}$ hence $g(z)=g_e(z)+\sum_{n\geq 1}\frac{g_n}{z^n}$ where $g_e$ is an entire function and the Laurent series converges absolutely for $z\ne 0$. Since  $g(z)={{O}\left(z^{-1}\right)}$ for $z\to\infty$ then $g_e=0$. Since $g(z)={{O}\left(z^{\al}\right)}$ for $z\to 0$ then $g(z)=\sum_{1\leq n\leq -\Re\al}\,{g_n}/{z^n}$. Since, more precisely, $g(z)=z^{\al}\mathcal{O}_0+{O}\left(1\right)$ for $z\to 0$ and $\al\not\in -\NN$ then $g(z)\equiv 0$, establishing uniqueness. \qed

\subsection{A Riemann-Hilbert problem on $[0,1]$}

Obviously, after a linear change of variable we can define the Hadamard singular part at any point: if, say, we change $t$ to $1-t$, then we can define for $f$ analytic at $1$
$$\nint_x^1(1-t)^\beta f(t)\,dt\, :=\, \ninto^{1-x}s^\beta\, f(1-s)\, ds$$
representing the analytic continuation of the usual integral, defined for $\beta\in\CC\setminus\ZZ_-$, and having the usual properties of integrals.

Let $\al, \beta\in\CC\setminus \ZZ_-$ (possibly with real parts less than $-1$). 

For functions $f$ analytic on $[0,1]$ define
\be\label{split3}\ninto^1 t^\al (1-t)^\beta f(t)\, dt\, :=\, \ninto^{\xi_0}  t^\al (1-t)^\beta f(t)\, dt\, \ \ \ \ \ \ \ \ \ \ \ \ \ \ \ee
$$\ \ \ \ \ \ \ \ \ \ \ \ \ \ +\, \int_{\xi_0}^{\xi_1}  t^\al (1-t)^\beta f(t)\, dt+ \nint_{\xi_1}^1t^\al (1-t)^\beta f(t)\, dt$$
where $0<\xi_0\leq \xi_1<1$ are any numbers so that $f$ is analytic on the disks $|z|<\xi_0$ and $|1-z|<1-\xi_1$, and the usual branches of $t^\al$ and $(1-t)^\beta$ are chosen. 

Formula (\ref{split3}) gives the analytic continuation in $\al$ and $\beta$ of the integral $\int_0^1t^\al (1-t)^\beta f(t)\, dt$ beyond the region $\Re\al>-1,\,\Re\beta>-1$.


{\em{Denote}} by $\mathcal{O}_1$ the germs of functions analytic at $1$.

 Let $p(x)$ be a polynomial.  
 Consider the following Riemann-Hilbert problem on $(0,1)$:

\z {{(i')}} $f$ is analytic on $\CC\setminus[0,1]$ and the following limits exist for $x\in(0,1)$:
$$f^+(x)=\lim_{z\to x;\Re z>0} f(z),\ \ \ \ \ f^-(x)=\lim_{z\to x;\Re z<0} f(z)$$

\z {{(ii')}} $\displaystyle{f^+(x)=f^-(x)+x^\al (1-x)^\beta p(x)}$ for $x\in(0,1)$.

\z {{(iii'$_\infty$)}} $f(z)=O(z^{-1})$ for $z\to\infty$,

\z {{(iii'$_0$)}} $f(z)=  z^\al\mathcal{O}_0+{O}(1)$ for $z\to 0$ and 

\z {{(iii'$_1$)}} $f(z)= (1-z)^\beta\mathcal{O}_1+{O}(1)$ for $z\to 1$.

\begin{Theorem}\label{ThRHj}
\z The Riemann-Hilbert problem {\em{(i')-(iii')}} has the unique solution $f=\mathcal{C}_{\al,\beta}[p]$ where
\be\label{notCab}
\mathcal{C}_{\al,\beta}[p]\, (z)=\frac{1}{2\pi i}\ninto^{1}\, \frac{t^\al(1-t)^\beta p(t)}{t-z}\, dt
\ee

\end{Theorem}

{\em{Proof of {Theorem}\,\ref{ThRHj}.}}

The proof is similar to that of {Theorem}\,\ref{ThRHsl}.

{\em{ Existence.}} Using a splitting (\ref{split3}) for the Hadamard finite part (\ref{notCab}) (where $\xi_0$ and $1-\xi_1$ are arbitrarily small) the finite part from $0$ to $\xi_0$ is analytic for $|z|>\xi_0$, and is  $O(z^{-1})$ for $z\to\infty$, the finite part from $\xi_1$ to $1$ is analytic for $|1-z|>1-\xi_1$, and is order $O(z^{-1})$ for $z\to\infty$, and by classical results, the middle integral is analytic in $z$ for $z\not\in[0,1]$, has the property {{(ii')}} for $\xi_0<x<\xi_1$, and satisfies {{(iii'$_\infty$)}}.

To establish the behavior of (\ref{notCab}) for $z\to 0$ note that $\mathcal{C}_{\al,\beta}[1]$ (for any $\al,\beta\not\in \ZZ_-$) satisfies the differential equation
\be\label{difxab}
\frac{d\phi}{dz}=\left(\frac{\al}{z}+\frac{\beta}{z-1}\right)\phi-\frac{a_0}{z}+\frac{a_1}{z-1}\ \ \ (a_{0,1}\ne 0\ {\mbox{constants}})
\ee
(see \S\ref{pfdifeqJ} for details) having the general solution $\phi(z)=\phi_{0;an}(z)+Cz^\al(z-1)^\beta $ with $\phi_{0;an}$ a function analytic at $z=0$. Therefore $\mathcal{C}_{\al,\beta}[1]= z^\al\mathcal{O}_0+{O}(1)$ for $z\to 0$. If $p(x)$ is a polynomial with $p(x)=\sum_{k=0}^{{\rm{deg}}p}p_{0,k}x^k$ then $\mathcal{C}_{\al,\beta}[p]=\sum_kp_{0,k}\mathcal{C}_{\al+k,\beta}[1]= z^\al\mathcal{O}_0+{O}(1)$ for $z\to 0$.

Similarly, the general solution of (\ref{difxab}) has the form $\phi(z)=\phi_{1;an}(z)+Cz^\al(z-1)^\beta $ with $\phi_{1;an}$ a function analytic at $z=1$, which implies that $\mathcal{C}_{\al,\beta}[p]=\sum_{k=0}^{{\rm{deg}}p}p_{1,k}\mathcal{C}_{\al,\beta+k}[1]= (1-z)^\beta\mathcal{O}_1+{O}(1)$ for $z\to 1$.

{\em{Uniqueness.}} Let $f$ be any function satisfying {{(i')}}-{{(iii')}}. Then $g=f-\mathcal{C}_{\al,\beta}[p]$
is analytic on $\CC\setminus[0,1]$ and $g^+(x)=g^-(x)$ for $x\in(0,1)$, therefore $g$ is analytic on $\CC\setminus \{0,1\}$. Take an integer $N>-\Re\al, -\Re\beta$. Then the function $\tilde{g}(z)=z^N(z-1)^Ng(z)$ is entire, since by (iii'$_0$) and (iii'$_1$) we have $\lim_{z\to 0}\tilde{g}(z)=0=\lim_{z\to 1}\tilde{g}(z)$. Therefore $g(z)$ has at most pole singularities at $z=0$ and $z=1$. Using again (iii'$_{0,1}$), since $\al,\beta\not\in\ZZ_-$ then $g$ is entire. Finally, by {{(iii'$_\infty$)}}, then $g=0$, proving uniqueness.
\qed

\section{Laguerre polynomials for general parameter}\label{Laguerre}

\subsection{Orthogonality of Laguerre polynomials for general parameter}\label{ortoL}

The definition used here for the Laguerre polynomials $L_n^{(\al)}$ is through their Rodrigues'  formula:
\be\label{RodriguesL}
L_n^{(\al)}(x)=w(x)^{-1}\frac{d^n}{dx^n}\left[ x^nw(x)\right]\ \ \ {\mbox{where\ }}w(x)=x^\al e^{-x}
\ee
which can be rewritten as a product of differential operators acting on the constant function $1$:
\be\label{LAprod}
L_n^{(\al)}(x)=\mathcal{A}_1\mathcal{A}_2\ldots\mathcal{A}_n 1,\ \ \ {\mbox{where\ }}\mathcal{A}_k=k+\al-x+x\partial_x
\ee
(this is not hard to show; see \cite{GenPol} for details). 

Using (\ref{LAprod}) it is clear that the leading monomial of $L_n^{(\al)}$ is $(-1)^nx^n$.

In the classical case, when $\Re\al>-1$, the Lagrange polynomials are orthogonal with respect to the bilinear functional
\be\label{BiLeg}
\calB_\al(f,g)=\int_0^\infty x^{\al}{\rm{e}}^{-x}f(x)g(x)\, dx
\ee

Analytic continuation in $\al$ of (\ref{BiLeg}) gives orthogonality in general:

\begin{Theorem}\label{genorto}
Let $\al\in\CC\setminus\ZZ_-$. Consider the bilinear functional on $\CC[x]$:
\be\label{nBiLeg}
\calB_\al(f,g)=\ninto^\infty x^{\al}{\rm{e}}^{-x}f(x)g(x)\, dx
\ee
We have
$$\calB_\al\left(L_n^{(\al)},L_k^{(\al)}\right)=0 \ {\mbox{for\ }}k\ne n\ \ \ \ \ {\mbox{and\ \ \ \ \ }}\calB_\al\left(L_n^{(\al)},L_n^{(\al)}\right)=n!\,\Gamma(\al+2)$$

\end{Theorem}

{{The proof}} of {Theorem}\,\ref{genorto} is immediate by analytic continuation, or directly, using the Rodrigues' formula (\ref{RodriguesL}) and integration by parts (Proposition\,\ref{FdThCal}).\qed

\begin{Corollary}\label{Coro1}
If $p_n$ is a polynomial of degree $n$ so that 
$$\calB_\al\left(p_n,x^k\right)=0\ {\mbox{for\ all\ }}k=0,1,\ldots n-1$$
 then $p_n$ is a scalar multiple of $L_n^{(\al)}$.
\end{Corollary}

\subsection{A Riemann-Hilbert problem for the Laguerre polynomials with general parameter}\label{RHLP}

\begin{Theorem}\label{ThRHL}

Let $\al\in\CC\setminus \ZZ_-$.

The following Riemann-Hilbert problem on $(0,+\infty)$ for $Y(z)\in\mathcal{M}_2(\CC)$:

\z {\em{(I)}} $Y$ is analytic on $\CC\setminus [0,+\infty)$ and the following limits exist for $x>0$:
$$Y^+(x)=\lim_{z\to x;\Re z>0} Y(z),\ \ \ \ \ Y^-(x)=\lim_{z\to x;\Re z<0} Y(z)$$

\z {\em{(II)}} $\displaystyle{Y^+(x)=Y^-(x)\left(\begin{array}{cc} 1& x^\al e^{-x}\\ 0 & 1\end{array}\right)}$ for $x>0$.

\z {\em{(III$_\infty$)}} $\displaystyle{Y(z)=\left( I+{O}\left(z^{-1}\right)\right)\, \left(\begin{array}{cc} z^n & 0\\ 0 & z^{-n}\end{array}\right)}$ for $z\to\infty$

\z {\em{(III$_0$)}} $\displaystyle{Y(z)= \left(\begin{array}{cc} {O}(1) & z^\al\mathcal{O}_0+{O}(1)\\  {O}(1) & z^\al\mathcal{O}_0+{O}(1) \end{array}\right)}$ for $z\to 0$

\z has the unique solution
$$Y=\left( \begin{array}{cc} \pi_n & \mathcal{C}_\al[\pi_n]\\ c_n\pi_{n-1} & c_n\mathcal{C}_\al[\pi_{n-1}] \end{array}\right)$$
where $\mathcal{C}_\al$ is the integral operator defined by (\ref{notC}),  $\pi_n$ are the monic Laguerre polynomials: $\pi_n=(-1)^nL_n^{(\al)}$ (with $L_n^{(\al)}$ defined by (\ref{RodriguesL})) and $c_n=-2\pi i\left[(n-1)!\Gamma(\al+n)\right]^{-1}$.

\end{Theorem}

{\em{Note}} that the asymptotic boundary conditions  {{(III)}} coincide with the classical ones when $\Re\al>0$ , or $-1<\al<0$. Note also that the expression for $Y$ is the analytic continuation in $\al$ of the one in the classical case.

{\em{Proof.}}
The proof of {Theorem}\,\ref{ThRHL} is similar to that in the classical case (see \cite{Deift_book}, \cite{Ismail}, \cite{KJ_lecture}). Details are provided here for completeness.

Since $Y^+_{11}=Y^-_{11}$ then $Y_{11}$ is analytic on $\CC\setminus \{ 0\}$ therefore $Y_{11}$ is given by a convergent McLaurin series $Y_{11}=\sum_{k\geq 0}a_kz^k+\sum_{k\geq 1}b_k/z^k$. Since $Y_{11}=O(1)$ ($z\to 0$) then all $b_k=0$ and since by {{(III$_\infty$)}} $Y_{11}=z^n(1+O(z^{-1}))$ ($z\to\infty$) then $Y_{11}=p_n$ a monic polynomial of degree $n$. 

Then $Y^+_{12}=Y^-_{12}+p_n(x)x^\al e^{-x}$ with $Y_{12}=O(z^{-n-1})$ for $z\to\infty$ and $Y_{12}=z^\al\mathcal{O}_0+O(1)$ for $z\to 0$, therefore, by Theorem\,\ref{ThRHsl}, we have $Y_{12}=\mathcal{C}_\al[p_n]$. To ensure better decay for $z\to\infty$ writing 
$$\frac{1}{t-z}=-\sum_{k=0}^{n-1}\frac{t^k}{z^{k+1}}+\frac{t^{n}}{z^{n}(t-z)}$$
 we have
$$\mathcal{C}_\al[p_n](z)=-\frac{1}{z^{k+1}}\sum_{k=0}^{n-1}\frac{1}{2\pi i}\ninto^{\infty}\, {t^k p_n(t)t^\al e^{-t}}\, dt$$
\be\label{ooeeqq}
+\frac{1}{z^{n}}\frac{1}{2\pi i}\ninto^{\infty}\, \frac{t^n p_n(t)t^\al e^{-t}}{t-z}\, dt\ee

The last Hadamard finite part in (\ref{ooeeqq}) is $O(z^{-1})$ by Theorem\,\ref{ThRHsl}. Therefore, the condition that  $Y_{12}=O(z^{-n-1})$ for $z\to\infty$ is equivalent to 
$$\ninto^{\infty}\, {t^k p_n(t)t^\al e^{-t}}\, dt=0\ {\mbox{for\ all\ }}k=0,1,\ldots,n-1$$
By {Corollary}\,\ref{Coro1} then $p_n$ is a multiple of $L_n^{(\al)}$, and since $p_n$ is monic then $p_n=\pi_n$ and $Y_{12}=\mathcal{C}_\al[\pi_n]$.

Next, we have $Y^+_{21}=Y^-_{21}$, with $Y_{21}=O(z^{n-1})$ for $z\to\infty$ and $Y_{21}=O(1)$ for $z\to 0$ therefore $Y_{21}$ is a polynomial of degree at most $n-1$: $Y_{21}=q_{n-1}$.

Finally, $Y^+_{22}=Y^-_{22}+q_{n-1}(x)x^\al e^{-x}$ with $Y_{22}=O(z^{-n})$ for $z\to\infty$ and $Y_{22}=z^\al\mathcal{O}_0+O(1)$ for $z\to 0$. By Theorem\,\ref{ThRHsl} then $Y_{22}=\mathcal{C}_\al[q_{n-1}]$. Since we must have $Y_{22}=z^{-n}(1+O(z^{-1}))$, then in the decomposition (\ref{ooeeqq}) of $\mathcal{C}_\al[q_{n-1}]$ we must have 
\be\label{situ1}
\ninto^{\infty}\, {t^k q_{n-1}(t)t^\al e^{-t}}\, dt=0\ {\mbox{for\ }}k=0,1,\ldots,n-2
\ee
and
\be\label{situ2} -\frac{1}{2\pi i}\,\ninto^{\infty}\, {t^{n-1} q_{n-1}(t)t^\al e^{-t}}\, dt=1
\ee
By {Corollary}\,\ref{Coro1} relations (\ref{situ1}) imply that $q_{n-1}=c_n\pi_{n-1}=(-1)^{n-1}c_nL_{n-1}^{(\al)}$ and then (\ref{situ2}) gives the stated value for $c_n$.\qed

\section{The Jacobi polynomials for general parameters}\label{Jacobi}

\subsection{Orthogonality of the Jacobi polynomials for general parameters}\label{ortoJ}

The Jacobi polynomials $P_n^{(\al,\beta)}$ are considered here on the interval $[0,1]$,  defined by the Rodrigues'  formula:
\be\label{RodrigJ}
P_n^{(\al,\beta)}(x)=w(x)^{-1}\frac{d^n}{dx^n}\left[ (x-x^2)^nw(x)\right]\ \ \ {\mbox{where\ }}
w(x)=x^\al (1-x)^\beta
\ee
which can be rewritten as a product of differential operators acting on the constant function $1$:
\begin{multline}\label{JAprod}
P_n^{(\al,\beta)}(x)=\mathcal{A}_1\ldots\mathcal{A}_n 1 \ \ {\mbox{where\ }}\\
\mathcal{A}_k=k+\al-x(2k+\al+\beta)+(x-x^2)\partial_x
\end{multline}
(see \cite{GenPol} for details). 

Using (\ref{JAprod}) it is clear that the leading monomial of $P_n^{(\al,\beta)}$ is found as $\prod_{k=1}^n(-(2k+\al+\beta)x-x^2\partial_x)$, which equals $C_nx^n$ with
\be\label{formCn}
C_n=(-1)^n\prod_{k=1}^n(k+n+\al+\beta)
\ee

Note that if $\al+\beta\in\{ -n-1,-n-2,-n-3,\ldots\}$ then the polynomial $P_n^{(\al,\beta)}$ has degree less than $n$.

In the classical setting, when $\Re\al, \Re\beta>-1$, the Jacobi polynomials are orthogonal with respect to the bilinear functional
\be\label{BiJa}
\calB_{\al,\beta}(f,g)=\int_0^1 x^\al (1-x)^\beta f(x)g(x)\, dx
\ee

Analytic continuation of (\ref{BiJa}) in the parameters $\al,\,\beta$ gives orthogonality in general:

\begin{Theorem}\label{genJorto}
Let $\al,\beta\in\CC\setminus\ZZ_-$ with $\al+\beta\not\in\{ -2,-3,-4\ldots\}$.

Consider the bilinear functional on $\CC[x]$:
\be\label{BiJeg}
\calB_{\al,\beta}(f,g)=\ninto^1 x^{\al}(1-x)^\beta f(x)g(x)\, dx
\ee
We have
$$\calB_{\al,\beta}\left(P_n^{(\al,\beta)},P_k^{(\al,\beta)}\right)=0 \ {\mbox{for\ }}k\ne n$$
\z and
$$\calB_{\al,\beta}\left(P_n^{(\al,\beta)},P_n^{(\al,\beta)}\right)=(-1)^n\,n!\,C_n\,\frac{\Gamma(\al+n+1)\Gamma(\beta+n+1)}{\Gamma(\al+\beta+n+2)}\, \ne\, 0$$
with $C_n$ given by (\ref{formCn}).

\end{Theorem}

The results of {Theorem}\,\ref{genorto} follow by analytic continuation in $\al,\beta$ of the corresponding results in the classical case. Alternatively, they can be immediately deduced using the Rodrigues' formula (\ref{RodriguesL}), integration by parts (Proposition\,\ref{FdThCal}) and (\ref{Bfun}).\qed

\begin{Corollary}\label{Coro2}
Let $\al,\beta\in\CC\setminus\ZZ_-$ with $\al+\beta\not\in\{ -2,-3,-4\ldots\}$.

If $p_n$ is a polynomial of degree $n$ so that 
$$\calB_{\al,\beta}\left(p_n,x^k\right)=0\ {\mbox{for\ all\ }}k=0,1,\ldots n-1$$
 then $p_n$ is a multiple of $P_n^{(\al,\beta)}$.
\end{Corollary}

\subsection{A Riemann-Hilbert problem for the Jacobi polynomials with general parameters}\label{RHJP}

\begin{Theorem}\label{ThRHJ}

Let $\al,\beta\in\CC\setminus\ZZ_-$ with $\al+\beta\not\in\{ -2,-3,-4\ldots\}$.

The Riemann-Hilbert problem on $\CC\setminus [0,1]$ for $Y(z)\in\mathcal{M}_2(\CC)$:

\z {\em{(I')}} $Y$ is analytic on $\CC\setminus[0,1]$ and the following limits exist for $x\in(0,1)$:
$$Y^+(x)=\lim_{z\to x;\Re z>0} Y(z),\ \ \ \ \ Y^-(x)=\lim_{z\to x;\Re z<0} Y(z)$$

\z {\em{(II')}} $\displaystyle{Y^+(x)=Y^-(x)\left(\begin{array}{cc} 1& x^\al (1-x)^\beta \\ 0 & 1\end{array}\right)}$ for $x\in (0,1)$.

\z {\em{(III'$_\infty$)}} $\displaystyle{Y(z)=\left( I+{O}\left(z^{-1}\right)\right)\, \left(\begin{array}{cc} z^n & 0\\ 0 & z^{-n}\end{array}\right)}$ for $z\to\infty$

\z {\em{(III'$_0$)}} $\displaystyle{Y(z)= \left(\begin{array}{cc} {O}(1) & z^\al\mathcal{O}_0+{O}(1)\\  {O}(1) & z^\al\mathcal{O}_0+{O}(1) \end{array}\right)}\ \ \ {\mbox{for\ }}z\to 0$

\z {\em{(III'$_1$)}} $\displaystyle{Y(z)= \left(\begin{array}{cc} {O}(1) & (z-1)^\beta\mathcal{O}_1+{O}(1)\\  {O}(1) & (z-1)^\beta\mathcal{O}_1+{O}(1) \end{array}\right)}\ \ \ {\mbox{for\ }}z\to 1$

\z has the unique solution
$$Y=\left( \begin{array}{cc} \pi_n & \mathcal{C}_{\al,\beta}[\pi_n]\\ c_n\pi_{n-1} & c_n\mathcal{C}_{\al,\beta}[\pi_{n-1}] \end{array}\right)$$
where $\mathcal{C}_{\al,\beta}$ is the integral operator (\ref{notCab}),  $\pi_n$ are the monic Jacobi polynomials: $\pi_n=C_n^{-1}P_n^{(\al,\beta)}$ and $c_n=(-1)^n2\pi i C_{n-1}\left[(n-1)!B(\al+n,\beta+n)\right]^{-1}$ (with $P_n^{(\al,\beta)}$ defined by (\ref{RodrigJ}), $C_n$ given by (\ref{formCn}), and $B$ given by (\ref{Bfun})).

\end{Theorem}

Note that the asymptotic boundary conditions  (III') are the classical ones if $0\ne\Re\al,\Re\beta>-1$.

 {\em{Proof.}}

The {{proof of {Theorem}\,\ref{ThRHJ}}} is similar to that in the classical case (see \cite{Deift_book}, \cite{Ismail}, \cite{KJ_lecture}). Details are provided here for completeness.

As in the proof of {Theorem}\,\ref{ThRHL} we obtain that $Y_{11}=p_n$ a monic polynomial of degree $n$. 

Then from $Y^+_{12}=Y^-_{12}+p_n(x)x^\al (1-x)^\beta$ with $Y_{12}=O(z^{-n-1})$ for $z\to\infty$ and $Y_{12}=z^\al\mathcal{O}_0+O(1)$ for $z\to 0$, and $Y_{12}=(z-1)^\beta\mathcal{O}_1+O(1)$ for $z\to 1$, by Theorem\,\ref{ThRHJ}, we have $Y_{12}=\mathcal{C}_{\al,\beta}[p_n]$. To ensure better decay for $z\to\infty$ we need, as in the proof of {Theorem}\,\ref{ThRHL}, that
\be\label{Jooeeqq}
\ninto^{1}\, \frac{t^k p_n(t)t^\al (1-t)^\beta}{t-z}\, dt\ \ \ {\mbox{for\ all\ }}k=0,1,\ldots,n-1
\ee

By {Corollary}\,\ref{Coro2} then $p_n$ is a multiple of $P_n^{(\al,\beta)}$ and since $p_n$ is monic then $p_n=\pi_n$ and $Y_{12}=\mathcal{C}_{\al,\beta}[\pi_n]$.

Next, we have $Y^+_{21}=Y^-_{21}$, with $Y_{21}=O(z^{n-1})$ for $z\to\infty$ and $Y_{21}=O(1)$ for $z\to 0$ therefore $Y_{21}$ is a polynomial of degree at most $n-1$: $Y_{21}=q_{n-1}$.

Finally, $Y^+_{22}=Y^-_{22}+q_{n-1}(x)x^\al (1-x)^\beta$ with $Y_{22}=O(z^{-n})$ for $z\to\infty$, $Y_{22}=z^\al\mathcal{O}_0+O(1)$ for $z\to 0$ and $Y_{22}=(z-1)^\beta\mathcal{O}_1+O(1)$ for $z\to 1$. By Theorem\,\ref{ThRHJ} then $Y_{22}=\mathcal{C}_{\al,\beta}[q_{n-1}]$. Since we must have $Y_{22}=z^{-n}(1+O(z^{-1}))$, then  
\be\label{situ1J}
\ninto^{1}\, {t^k q_{n-1}(t)t^\al (1-t)^\beta}\, dt=0\ {\mbox{for\ all\ }}k=0,1,\ldots,n-2
\ee
and
\be\label{situ2J}
-\frac{1}{2\pi i}\, \ninto^{1}\, {t^{n-1} q_{n-1}(t)t^\al (1-t)^\beta}\, dt=1
\ee
By {Corollary}\,\ref{Coro2} relations (\ref{situ1J}) imply that $q_{n-1}=c_n\pi_{n-1}=c_nC_{n-1}^{-1}L_{n-1}^{(\al)}$, then (\ref{situ2J}) gives the stated value for $c_n$.\qed



\section{Appendix}

\subsection{Proof of Proposition\,\ref{addi}.}\label{pfaddi}
The equality clearly holds for $f(t)=t^n,\,n\in\NN$ and follows for all $f$ by linearity, density of polynomials in $\calB(D_0)$ and continuity of the operator $\calJ_{\al}$.\qed

\subsection{Proof of Proposition\,\ref{FdThCal}.}\label{pfFdThCal}

Note that $\frac{d}{dt}\left[ t^\al f(t)\right]=t^{\al-1}[\al f(t)+t f'(t)]$, so the definition (\ref{def_int}) applies.

A simple calculation using (\ref{defJ}) and (\ref{def_int})  yields
\begin{multline}\nonumber
{ \ninto}^{x}\, \frac{d}{dt}\left[ t^\al f(t)\right]=\ninto^{x}\, \frac{d}{dt}\left[ t^\al\sum_{n=0}^\infty f_nt^{n}\right]\, dt={ \nint}_{0}^{x}\, t^{\al-1}\sum_{n=0}^\infty (n+\al)f_nt^{n}\\
=\sum_{n=0}^\infty f_nx^n
\end{multline}
\qed

\subsection{$\mathcal{C}_\beta[1]$ satisfies the equation (\ref{difxbeta})}\label{pfdifeq}

Denote
$$f(z)={2\pi i}\,\mathcal{C}_\beta[1](z)=\ninto^{\infty}\, \frac{t^\beta e^{-t}}{t-z}\, dt$$
Differentiation and integration by parts give
$$f'(z)=\ninto^{\infty}\, \frac{t^\beta e^{-t}}{(t-z)^2}\, dt=\ninto^{\infty}\, \frac{t^\beta(\beta/t-1)e^{-t}}{t-z}\, dt$$
$$=\ninto^{\infty}\, \frac{\beta}{z}\left(\frac{1}{t-z}-\frac{1}{t}\right){t^\beta e^{-t}}\, dt -\ninto^{\infty}\, \frac{t^\beta e^{-t}}{t-z}\, dt$$
which gives (\ref{difxbeta}).

\subsection{$\mathcal{C}_{\al,\beta}[1]$ satisfies the equation (\ref{difxab})}\label{pfdifeqJ}

Denote
$$f(z)={2\pi i}\,\mathcal{C}_{\al,\beta}[1]=\ninto^{1}\, \frac{t^\al(1-t)^\beta }{t-z}\, dt$$
Differentiation and integration by parts give
$$f'(z)=\ninto^{1}\, \frac{t^\al(1-t)^\beta }{(t-z)^2}\, dt =\ninto^{1}\, \frac{t^\al(1-t)^\beta }{t-z}\left[\frac{\al}{t}-\frac{\beta}{1-t}\right]\, dt$$
$$=\ninto^{1}\, {t^\al(1-t)^\beta }\left[\frac{\al}{z}\left(\frac{1}{t-z}-\frac{1}{t}\right)+\frac{\beta}{z-1}\left(\frac{1}{t-z}+\frac{1}{1-t}\right)\right]\, dt$$
which gives (\ref{difxab}) with $a_0=\frac{\al}{2\pi i}\, B(\al,\beta+1)$ and $a_1=\frac{\beta}{2\pi i}\, B(\al+1,\beta)$ where 
\be\label{Bfun}
B(r,s)=\ninto^{1}\, t^{r-1}(1-t)^{s-1}\, dt=\frac{\Gamma(r)\Gamma(s)}{\Gamma(r+s)}
\ee
(the last equality in (\ref{Bfun}) is valid by analytic continuation beyond the region $\Re r,\Re s>0$).

\

\end{document}